\documentclass[12pt]{article}

\usepackage{amsmath,amssymb,amsfonts,amsthm}
\usepackage{tikz,color,pgf,graphicx,url}

\newtheorem{theorem}{Theorem}[section]
\newtheorem{corollary}[theorem]{Corollary}
\newtheorem{lemma}[theorem]{Lemma}
\newtheorem{proposition}[theorem]{Proposition}

\renewcommand{\proof}{\noindent{\bf Proof.\ }}
\renewcommand{\qed}{\hfill $\square$ \bigskip}

\newcommand{\ggcd}{\gamma_{\rm cg}}
\newcommand{\ggcs}{\gamma_{\rm cg}'}
\newcommand{\ggd}{\gamma_{\rm g}}
\newcommand{\ccg}{\! \rightarrow \!}
\newcommand{\cG}{{\cal G}}

\tikzstyle{every node}=[circle, draw, fill=black!10, inner sep=0pt, minimum width=4pt]

\newcommand{\fast}{{\sc Fast}}
\newcommand{\slow}{{\sc Slow}}

\begin{document}

\title{Predominating a vertex in the connected domination game}

\author{Csilla Bujt\'as$^{a,d}$\thanks{Email: \texttt{csilla.bujtas@fmf.uni-lj.si}} 
\and Vesna Ir\v si\v c$^{a,b,e}$\thanks{Email: \texttt{vesna.irsic@fmf.uni-lj.si}}
\and Sandi Klav\v zar $^{a,b,c}$\thanks{Email: \texttt{sandi.klavzar@fmf.uni-lj.si}}
}
\maketitle

\begin{center}
$^a$ Faculty of Mathematics and Physics, University of Ljubljana, Slovenia\\
\medskip

$^b$ Institute of Mathematics, Physics and Mechanics, Ljubljana, Slovenia\\
\medskip

$^c$ Faculty of Natural Sciences and Mathematics, University of Maribor, Slovenia\\
\medskip

$^d$ Faculty of Information Technology, University of Pannonia, Hungary\\
\medskip

$^e$ Department of Mathematics, Simon Fraser University, Burnaby, BC, Canada\\
\medskip

\end{center}

\begin{abstract}
The connected domination game is played just as the domination game, with an additional requirement that at each stage of the game the vertices played induce a connected subgraph. The number of moves in a D-game (an S-game, resp.) on a graph $G$ when both players play optimally is denoted by $\gamma_{\rm cg}(G)$ ($\gamma_{\rm cg}'(G)$, resp.). Connected Game Continuation Principle is established as a substitute for the classical Continuation Principle which does not hold for the connected domination game. Let $G|x$ denote the graph $G$ together with a declaration that the vertex $x$ is already dominated. The first main result asserts that if $G$ is a graph with $\gamma_{\rm cg}(G) \geq 3$ and $x \in V(G)$, then $\gamma_{\rm cg}(G|x) \leq 2 \gamma_{\rm cg}(G) - 3$ and the bound is sharp. The second main theorem states that if $G$ is a graph with $n(G) \geq 2$ and $x \in V(G)$, then $\gamma_{\rm cg}(G|x) \geq \left \lceil \frac12 \gamma_{\rm cg}(G) \right \rceil$ and the bound is sharp. Graphs $G$ and their vertices $x$ for which $\gamma_{\rm cg}'(G|x) = \infty$ holds are also characterized. 
\end{abstract}


\noindent
{\bf AMS Subj.\ Class.\ (2020)}: 05C57, 05C69

\subsection*{Declarations}

\subsubsection*{Funding} 
We acknowledge the financial support from the Slovenian Research Agency (research core funding No.\ P1-0297 and projects J1-9109, J1-1693, N1-0095, N1-0108). 

\subsubsection*{Conflicts of interest} 
Not applicable.

\subsubsection*{Availability of data and material}
Not applicable.

\subsubsection*{Code availability} 
The code used was written in python and is not publicly available.

\section{Introduction}

The {\em domination game}~\cite{bresar-2010} is played on a graph $G$ by Dominator and Staller who take turns, each time selecting a vertex which dominates at least one vertex that has not yet been dominated by the vertices already played. The game is over when no move is possible. The goal of Dominator is to select as few vertices as possible, Staller's goal is just the opposite. Assuming that both players are playing optimally, the number of vertices selected by the end of the game is a graph invariant. If Dominator has the first move, then the invariant is called the {\em game domination number} of $G$, denoted by $\ggd(G)$, otherwise (if Staller starts the game) it is denoted by $\ggd'(G)$. The total domination game~\cite{henning-2015} which is played on an isolate-free graph follows the same rules, except that each newly played vertex must totally dominate at least one new vertex. The corresponding invariants are denoted by $\gamma_{\rm tg}(G)$ and $\gamma_{\rm tg}'(G)$. For some recent results on the (total) domination game see~\cite{bujtas-2021b, charoensitthichai-2020, jiang-2019}, for a variety of the classical domination game see~\cite{bresar-2019}, for the fractional domination game~\cite{fractional-2019}, for Maker-Breaker domination games~\cite{duchene-2020, gledel-2020}, and for a state of the art on domination games till the early 2021 the book~\cite{book-2021}. 
 
In this paper we are intrigued by the \emph{connected domination game} which was introduced by Borowiecki, Fiedorowicz, and Sidorowicz~\cite{borowiecki+2019connected}. The game is played on a connected graph and the rules of the game are the same as for the
 domination game, except that a move is legal if the selected vertex not only dominates a vertex which is not yet dominated by previous moves, but is also adjacent to at least one already played vertex. Note that the latter requirement is equivalent to the fact that at each stage of the game, the set of vertices played induces a connected subgraph. If Dominator starts the connected domination game, then we call it a \emph{D-game}, while if Staller is the first to select a vertex, then we speak of an \emph{S-game}. If both players play optimally, the number of moves in a D-game is the \emph{game connected domination number} $\ggcd(G)$, and the number of moves in an S-game is the \emph{Staller-start game connected domination number} $\ggcs(G)$. 

If $S \subseteq V(G)$, then a \emph{partially dominated graph} $G|S$ is a graph together with a declaration that the vertices from $S$ are already dominated, that is, the vertices from $S$ need not be dominated during the course of the (total/connected) domination game. If $S = \{u\}$, then the notation is simplified to $G|u$. We say that a connected domination game on $G$ or $G|u$ is \emph{optimal} if it is a sequence of moves such that both players play according to their optimal strategies. We will use the convention to denote the sequence of moves in a D-game by $d_1, s_1, \ldots$, and by $s_1', d_1', \ldots$ in an S-game. 

In~\cite[Theorem 3]{bujtas-2015} it was proved that if $u$ is a vertex of a graph $G$, then $\ggd(G)\le \ggd(G|u) + 2$. On the other hand, the Continuation Principle~\cite{kinnersley-2013} implies that $\ggd(G|u) \le \ggd(G)$ holds. We thus have:
\begin{equation}
\label{eq:usual-game}
\ggd(G) - 2 \le \ggd(G|u) \le \ggd(G)\,.
\end{equation}
Similarly, for the total domination game it was proved in~\cite[Lemma 2.1]{irsic-2019} that if $u$ is a vertex of a graph $G$ that contains no isolated vertices, then $\gamma_{\rm tg}(G)  \le \gamma_{\rm tg}(G|u) + 2$. Since the Continuation Principle holds for the total domination game as well~\cite{henning-2015}, for a graph $G$ without isolated vertices we have
\begin{equation}
\label{eq:total-game}
\gamma_{\rm tg}(G) - 2 \le \gamma_{\rm tg}(G|u) \le \gamma_{\rm tg}(G)\,.
\end{equation}
All the bounds in~\eqref{eq:usual-game} and~\eqref{eq:total-game} are sharp~\cite{bujtas-2015, irsic-2019}. For more information on the game (total) domination number of graphs with one vertex predominated see~\cite{dorbec-2019, henning-2018a, henning-2018b, xu-2018}.

Answering~\cite[Problem 6.2]{irsic2019+connected}, our main results read as follows. 

\begin{theorem}
	\label{thm:G|x_upper}
	If $G$ is a connected graph with $\ggcd(G) \geq 3$ and $x \in V(G)$, then $$\ggcd(G|x) \leq 2 \ggcd(G) - 3\,.$$
	Moreover, the bound is sharp. 
\end{theorem}

\begin{theorem}
	\label{thm:G|x_lower}
	If $G$ is a connected graph on at least two vertices and $x \in V(G)$, then $$\ggcd(G|x) \geq \left \lceil \frac12 \ggcd(G) \right \rceil.$$
	Moreover, the bound is sharp. 
\end{theorem}

Comparing these two theorems with~\eqref{eq:usual-game} and~\eqref{eq:total-game} reveals that the connected domination game is very different from the (total) domination game. This difference in particular follows from the fact that the Continuation Principle in the usual sense does not hold for the connected domination game. 

The paper is organized as follows. In the next section we give additional definitions, recall some results, and provide a new family of sharpness examples for the earlier established upper bound $\ggcs(G) \leq 2 \ggcd(G)$. We refer to this construction in our further developments. In Section~\ref{sec:continuation} we prove, what we call, Connected Game Continuation Principle. We consider it as a substitute for the usual Continuation Pronciple and apply it in the continuation of the paper. In Sections~\ref{sec:proof-of-first-theorem} and~\ref{sec:proof-of-second-theorem}, Theorems~\ref{thm:G|x_upper} and~\ref{thm:G|x_lower} are proved, respectively. In Section~\ref{sec:S-game} we turn our attention to the S-game and characterize the graphs $G$ and its vertices $x$ for which $\ggcs(G|x) = \infty$ holds. In particular, if $x$ is a vertex of a tree $T$, then $\ggcs(T|x) = \infty$ if and only if $x$ is not a leaf and has a neighbor of degree $2$.

\section{Preliminaries}
\label{sec:prelim}

Let $G$ be a graph. If $S \subseteq V(G)$, then the subgraph induced by $S$ is denoted by $G[S]$. For a vertex $v \in V(G)$, the \emph{(open) neighborhood} $N(v)$ is the set of neighbors of $v$, and the \emph{closed neighborhood} $N[v] = N(v) \cup \{v\}$. If $S \subseteq V(G)$, then $N[S] = \bigcup_{v\in S} N[v]$. A vertex $v \in V(G)$ \emph{dominates} itself and its neighbors. A subset of vertices $D \subseteq V(G)$ is a \emph{dominating set} of $G$ if it dominates all vertices of $G$, i.e.\ $N[D] = V(G)$. This means that every vertex from $V(G) \setminus D$ has a neighbor in $D$. The minimum cardinality of a dominating set of $G$ is the \emph{domination number} $\gamma(G)$ of $G$. Similarly, a vertex $v \in V(G)$ \emph{totally dominates} its neighbors, but not itself. A \emph{total dominating set} of $G$ is a subset $D \subseteq V(G)$ if every vertex from $V(G)$ has a neighbor in $D$. Minimum cardinality of a total dominating set of an isolate-free graph $G$ is the \emph{total domination number} $\gamma_{\rm t}(G)$. A \emph{connected dominating set} $D$ of $G$ is a dominating set such that $G[D]$ is connected. Minimum cardinality of such a set in a connected graph $G$ is the \emph{connected domination number} $\gamma_{\rm c}(G)$. For a positive integer $n$ we use notation $[n] = \{1, \ldots, n\}$.

If $x$ is a vertex of a connected graph $G$, then $\ggcd(G|x)$ is well-defined. Indeed, Dominator can play $x$ as his first move, and the rest of the game is then a usual connected domination game which always finishes after a finite number of moves. However, this may not be the case for the S-game. As a simple example consider the path $P_4$ on vertices $v_1, v_2, v_3, v_4$ with natural adjacency relation, and consider the S-game played on $P_4|v_2$. Then, after the first move $s_1' = v_4$, Dominator has no legal vertex to play. This means that the game cannot be finished as $v_1$ remains  undominated. Because of this phenomenon we will write $\ggcs(G|x) = \infty$ if Staller has a strategy in the S-game played on $G|x$ such that at some point of the game no legal moves are available, but not all vertices are already dominated. 

The following basic property of the game connected domination number will be useful.

\begin{theorem}[{\cite[Theorem 1]{borowiecki+2019connected}}]
	\label{thm:c-bound}
	If $G$ is a graph, then $$\gamma_c(G) \leq \ggcd(G) \leq 2 \gamma_c(G) -1.$$
\end{theorem}

In our later arguments, the graphs $G$ for which $\ggcd(G) = \gamma_{\rm c}(G)$ holds will be important. In this respect we mention that forests $F$ for which $\ggd(F) = \gamma(F)$ holds were characterized in~\cite{nadjafi-2016}, while trees with equal total domination and game total domination numbers were described in~\cite{henning-2017}. For a more general framework in this direction, see~\cite{bujtas-2021}. 

As stated, Theorem~\ref{thm:c-bound} was proved in~\cite{borowiecki+2019connected}, but a more detailed proof of the upper bound was later presented in~\cite{bujtas+2019connected}. Its proof reveals the following fact that we state here for later usage.

\begin{lemma}
\label{lem:Dom-strategy-inD-game}
Let $G$ be a connected graph and let $S$ be a connected dominating set of $G$. For every $v \in S$, Dominator has a strategy to start a D-game on $G$ by playing $v$, and playing only vertices from $S$ during the game.
\end{lemma}

Note that we can use Lemma~\ref{lem:Dom-strategy-inD-game} for a (connected) subgraph of $G$ as well. The proofs of the following lemmas are analogous to that of Lemma~\ref{lem:Dom-strategy-inD-game} and hence we omit them here.

\begin{lemma}
	\label{lem:1}
	If $D$ is a connected dominating set of $G|x$, then Dominator has a strategy to play only vertices from $D$ during the game. Furthermore, the game ends after at most $2 |D| - 1$ moves.
\end{lemma}

\begin{lemma}
	\label{lem:2}
	Let $k \geq 0$ be even and let $D$ be the set of the first $k$ moves of a connected domination game. If $D' \subseteq V(G) \setminus D$ is a set such that $D \cup D'$ is a connected dominating set, then Dominator has a strategy to play only vertices from $D'$ during the remaining part of the game. Furthermore, the game ends after at most $|D| + 2 |D'| - 1$ moves.
\end{lemma}

We next recall the following result and demonstrate its sharpness.

 \begin{theorem}[{\cite[Theorem 3.2]{irsic2019+connected}}]
	\label{thm:staller-start}
	If $G$ is a graph, then $$\ggcs(G) \leq 2 \ggcd(G).$$
\end{theorem}

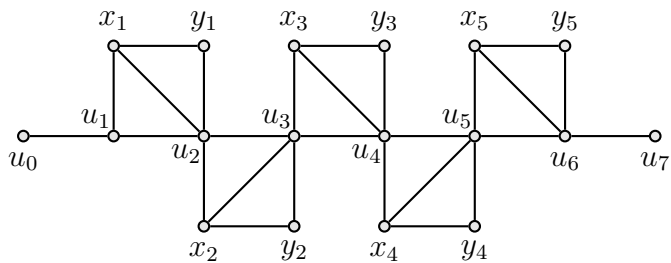
\begin{figure}[!ht]
	\begin{center}
		\begin{tikzpicture}[thick, scale=1.2pt]
		
		\node[label=below: {$u_0$}] (u0) at (0,0) {};
		\node[label=135: {$u_1$}] (u1) at (1,0) {};
		\node[label=-135: {$u_2$}] (u2) at (2,0) {};
		\node[label=135: {$u_3$}] (u3) at (3,0) {};
		\node[label=-135: {$u_4$}] (u4) at (4,0) {};
		\node[label=135: {$u_5$}] (u5) at (5,0) {};
		\node[label=below: {$u_6$}] (u6) at (6,0) {};
		\node[label=below: {$u_7$}] (u7) at (7,0) {};
		
		\node[label=above: {$x_1$}] (x1) at (1,1) {};
		\node[label=below: {$x_2$}] (x2) at (2,-1) {};
		\node[label=above: {$x_3$}] (x3) at (3,1) {};
		\node[label=below: {$x_4$}] (x4) at (4,-1) {};
		\node[label=above: {$x_5$}] (x5) at (5,1) {};
		
		\node[label=above: {$y_1$}] (y1) at (2,1) {};
		\node[label=below: {$y_2$}] (y2) at (3,-1) {};
		\node[label=above: {$y_3$}] (y3) at (4,1) {};
		\node[label=below: {$y_4$}] (y4) at (5,-1) {};
		\node[label=above: {$y_5$}] (y5) at (6,1) {};

		\draw (u0) -- (u1) -- (u2) -- (u3) -- (u4) -- (u5) -- (u6);
		
		\draw (u1) -- (x1) -- (y1);
		\path (u2) edge (x1);
		\path (u2) edge (y1);
		
		\draw (u2) -- (x2) -- (y2);
		\path (u3) edge (x2);
		\path (u3) edge (y2);
		
		\draw (u3) -- (x3) -- (y3);
		\path (u4) edge (x3);
		\path (u4) edge (y3);
		
		\draw (u4) -- (x4) -- (y4);
		\path (u5) edge (x4);
		\path (u5) edge (y4);
		
		\draw (u5) -- (x5) -- (y5);
		\path (u6) edge (x5);
		\path (u6) edge (y5);
		
		\path (u6) edge (u7);
		
		\end{tikzpicture}
		\caption{The graph $H_6$.}
		\label{fig:primer-H6}
	\end{center}
\end{figure}

Let $H_n$, $n \geq 2$, be a graph with vertices $V(G_n) = \{ u_0, \ldots, u_{n+1}\} \cup \{ x_1, \ldots, x_{n-1}\} \cup \{y_1,  \ldots, y_{n-1} \}$ and edges $u_i u_{i+1}$ for $i \in \{0, \ldots, n\}$, $u_i x_i$, $x_i y_i$, $y_i u_{i+1}$, and $u_{i+1}  x_i$ for $i \in [n-1]$. See Fig.~\ref{fig:primer-H6} for $H_6$. The family $H_n$ is actually obtained by a simplification of a family $G_n$ from~\cite{irsic2019+connected}. This simplification was proposed by West~\cite{westBled}.

We recall the strategy of Dominator from~\cite[Lemma 3.3]{irsic2019+connected}.  His strategy is to play $d_1 = u_n$ which makes all the remaining moves unique and the game finishes in $n$ moves. We call this strategy \fast. Note that exactly vertices $u_n, \ldots, u_1$ are played during the game. Together with the fact that $\gamma_{\rm c}(H_n) = n$, and hence $\ggcd(H_n) \geq n$ by Theorem~\ref{thm:c-bound}, we get:

\begin{lemma}
	\label{lem:Hn-Dom}
	If $n \geq 2$, then $\ggcd(H_n) = n$.
\end{lemma}

However, to determine $\ggcs(H_n)$, the method from the proof of \cite[Lemma 3.4]{irsic2019+connected} is not helpful, thus we use a different approach. 

\begin{lemma}
	\label{lem:Hn-St}
	If $n \geq 2$, then $\ggcs(H_n) = 2n$.
\end{lemma}

\proof
	It follows from Theorem~\ref{thm:staller-start} and Lemma~\ref{lem:Hn-Dom} that $\ggcs(H_n) \leq 2 n$. To prove the reverse inequality, we consider the following strategy of Staller. 
	
	She starts the game by playing $s_1' = u_0$. Dominator's only legal reply is $d_1' = u_1$. Now Staller can play $s_2' = x_1$, which forces Dominator to select $d_2' = u_2$. Similarly, for $k \in \{ 3, \ldots, n-1 \}$, Staller can play $s_k' = x_{k-1}$, which leaves only one possible reply for Dominator, $d_k' = u_k$. After $2(n-1)$ moves, all vertices except $y_{n-1}$ and $u_{n+1}$ are dominated. Next Staller can play $s_n' = x_{n-1}$, which leaves Dominator finishing the game by playing $d_n' = u_n$. This strategy of Staller ensures that $\ggcs(H_n) \geq 2n$.
\qed

We call the strategy of Staller described in the proof of Lemma~\ref{lem:Hn-St} \slow. In short, her strategy is to start on $u_0$ and to play vertices $x_1, \ldots, x_{n-1}$ whenever she can. She is able to force $n-1$ additional moves on $V(H_n) \setminus \{ u_0, \ldots, u_{n+1} \}$. Thus, apart from the vertices $u_1, \ldots, u_n$, exactly $n$ additional moves are played (counting the move $u_0$ as well).

Finally, we describe the connected domination game with Chooser. Its rules are the same as in the connected domination game, except that there is another player, Chooser, who can make zero, one or more moves after any move of Dominator or Staller. The only rule for his move to be legal is that the set of played vertices is still connected after his move. The following holds: 
\begin{lemma}[{Chooser Lemma~\cite{borowiecki+2019connected}}]
	\label{lema:chooser}
	Consider the connected domination game with Chooser on a graph $G$. Suppose that in the game Chooser plays $k$ vertices, and that both Dominator and Staller play optimally. Then at the end of the game the number of played vertices is at most $\ggcd(G) + k$ and at least $\ggcd(G) - k$.
\end{lemma}

\section{Connected Game Continuation Principle}
\label{sec:continuation}

Although the usual form of the Continuation Principle does not hold for the connected domination game, in this section we establish a variation of it under the condition that a game which has been already started must be continued. First, we state a lemma that can be proved with parallel arguments as Lemma~\ref{lem:Dom-strategy-inD-game}. 
\begin{lemma}
	\label{lem:legal sequence}
	Let $G$ be a  graph and let $S \subseteq V(G)$ such that $G[S]$ is connected. For every vertex $v \in S$, there is a sequence $(v_1=v, v_2, \dots, v_k)$ such that $v_i \in S$ for all $i \in [k]$, $N[\{v_1, \dots ,v_k\}]=N[S]$, and  $(v_1, \dots, v_k)$ is a legal sequence of the first $k$ moves in a  D-game on $G$.  
\end{lemma}
Under the conditions of Lemma~\ref{lem:legal sequence} but without specifying $v$,  we say that $S'=\{v_1, \dots ,v_k\}$ is a \emph{legal subset} of $S$ in $G$. Clearly, every set $S$ that induces a connected subgraph in $G$ has at least one legal subset. 
\medskip

Now, we introduce two invariants related to the optimal continuation of a connected domination game. Let $G$ be a connected graph and $D \subseteq V(G)$ such that $D\neq \emptyset$ and $G[D]$ is connected. Let $D'$ be a legal subset of $D$. The \emph{connected game continuation number} $\gamma_{\rm cg}(D \ccg G)$  is the number of moves needed to finish a connected domination game where the first $|D'|$ moves are the elements of $D'$ in any legal order, both players play optimally, and Dominator plays the $(|D'|+1)^{\rm st}$ move in the game. The invariant $\gamma_{\rm cg}'(D \ccg G)$ is defined analogously, if Staller selects the $(|D'|+1)^{\rm st}$ vertex in the game.  By definition,  if $D'$ and $D''$ are two different legal subsets of $D$ in $G$, then $N[D']=N[D'']$ and 
$$N[D']\setminus \{v\in V(G): N[v] \subseteq N[D']\}=
N[D'']\setminus \{v\in V(G): N[v] \subseteq N[D'']\}.$$
That is, the set of dominated vertices and that of playable vertices remain the same independently of the specification of the legal subset. This shows that the definitions of $\gamma_{\rm cg}(D \ccg G)$ and $\gamma_{\rm cg}'(D \ccg G)$ are sound.

Now the announced principle reads as follows.

\begin{lemma}[Connected Game Continuation Principle, CGCP]
	\label{lem:CGCP}
	Let $G$ be a connected graph, and let $D$ and $C$ be two nonempty subsets of $V(G)$ such that both $G[D]$ and $G[D \cup C]$ are connected. Then, $$\gamma_{\rm cg}((D\cup C) \ccg G) \le \gamma_{\rm cg}(D \ccg G) \quad \text{and} \quad \gamma_{\rm cg}'((D\cup C) \ccg G) \le \gamma_{\rm cg}'(D \ccg G).$$
\end{lemma}

The proof of Lemma~\ref{lem:CGCP} is basically the same as the one presented for the Chooser Lemma in~\cite{borowiecki+2019connected}. To be self contained, we present the main idea of the proof.

\noindent{\bf Sketch of the proof of Lemma~\ref{lem:CGCP}.\ }
	Dominator's strategy when the game is played on $G|N[D \cup C]$ is to imagine the game on $G|N[D]$, select his optimal move there, and play a corresponding move in the real game. Dominator plays in such a way that he maintains the following invariant: every dominated vertex in the imagined game is also dominated in the real game, and the number of played vertices in the real game is less than or equal to the number of played vertices in the imagined game. Proving that Dominator can always play such that he satisfies the above invariant is similar as the proof of~\cite[Chooser Lemma]{borowiecki+2019connected}.
\qed


As a consequence of the CGCP,  we may compare two possible legal moves $v$ and $v'$ during a connected domination game. Let $D'$ be the set of moves played until this point in the game. Assuming that $v'$ dominates all vertices which would be newly dominated by the move $v$, we can apply the CGCP with  $D=D' \cup \{v\}$ and $C=\{v'\}$ and get that 
$$\gamma_{\rm cg}((D'\cup  \{v,v'\}) \ccg G) \le \gamma_{\rm cg}((D' \cup \{v\}) \ccg G)$$
where, by our condition, $\gamma_{\rm cg}((D'\cup  \{v,v'\}) \ccg G)$ clearly corresponds to $\gamma_{\rm cg}((D'\cup  \{v'\}) \ccg G)$. This proves that playing $v'$ in the next turn is at least as advantageous for Dominator as playing $v$. If it is Staller's turn, playing $v$ is at least as advantageous for her as playing $v'$. We formulate this fact in the next result.

\begin{corollary}
	\label{cor:adv}
	Let $D$ be the set of the first $k \ge 0$ moves of a connected domination game and suppose that $v$ and $v'$ are two possible legal moves in the next turn such that $N[v] \cap (V(G)\setminus N[D]) \subseteq N[v'] \cap (V(G)\setminus N[D])$ holds.
	Then, 
	\begin{align*}
	\gamma_{\rm cg}((D\cup \{v'\}) \ccg G) &\le \gamma_{\rm cg}((D\cup \{v\}) \ccg G) \enskip \text{and}\\
	\gamma_{\rm cg}'((D\cup \{v'\}) \ccg G) &\le \gamma_{\rm cg}'((D\cup \{v\}) \ccg G).
	\end{align*}
\end{corollary}
In particular, as Corollary~\ref{cor:adv} immediately implies, if $N[v] \subseteq N[v']$, then Dominator may always play $v'$ instead of $v$ and Staller may always play $v$ instead of $v'$ during an optimal connected domination game on $G$.

\section{Proof of Theorem~\ref{thm:G|x_upper}}
\label{sec:proof-of-first-theorem}

Recall the statement of Theorem~\ref{thm:G|x_upper}: if $G$ is a connected graph with $\ggcd(G) \geq 3$ and $x \in V(G)$, then \begin{equation}
\label{eq:upper}
\ggcd(G|x) \leq 2 \ggcd(G) - 3\,,
\end{equation}
and the bound is sharp. Note that the condition $\ggcd(G) \geq 3$ is necessary as~\eqref{eq:upper} cannot hold if $\ggcd(G) = 1$, and there are plenty of graphs $G$ with $\ggcd(G) = 2$ and $\ggcd(G|x) > 1$ (for example a graph obtained from two copies of $K_n$ ($n \ge 3$) with exactly one edge between them). Before proving the theorem we add that a slightly weaker bound 
\begin{equation*}
\ggcd(G|x) \leq 2 \ggcd(G) - 1
\end{equation*}
is significantly easier to obtain---see~\eqref{eq:weak-bounds} below.

\subsection{Proof of~\eqref{eq:upper}}

In the proof, we consider a fixed connected graph $G$ with a fixed vertex $x \in V(G)$ and distinguish six cases concerning $G$ and $x$. For each case, we describe a strategy for Dominator which ensures that the connected domination game on $G|x$ finishes with at most $2\ggcd(G)-3$ played vertices. The moves in the connected domination game on $G|x$ will be denoted by $d_1^*,s_1^*, \dots$
When Case $i$ is discussed, for each $i \ge 2$, we suppose that none of the conditions of the previous cases can be applied.

For the graph $G$, we define $\cG$ as the set of all optimal D-games on $G$. A game $P \in \cG$ is represented by the sequence $(d_1,s_1, \dots)$ of moves. Each $P\in \cG$ consists of exactly $\ggcd(G)$ entries, the set of the first $i$ of them is denoted by $\widehat{P}_i$,  for all $i \in [\ggcd(G)]$. We will write $\widehat{P}$ instead of $\widehat{P}_{\ggcd(G)}$. Let us set $k=\lceil \frac{\ggcd(G)}{2} \rceil$. If $\ggcd(G)=2k$, then $P=(d_1, s_1, \dots ,d_k, s_k)$, whilst in case of $\ggcd(G)=2k-1$, we have $P=(d_1, s_1, \dots ,d_k)$. 

\begin{description}
	\item[Case 1:] $\gamma_{\rm c}(G|x) \le \gamma_{\rm c}(G)-1 $ or $\gamma_{\rm c}(G) \le \ggcd(G)-1$.
	
	First, we show that the inequality chain 
	\begin{equation} \label{eq:weak-bounds}
	\ggcd(G|x) \le 2 \gamma_{\rm c}(G|x) - 1 \le 2 \gamma_{\rm c}(G) - 1 \le 2 \ggcd(G)-1.
	\end{equation}
	always holds, independently of the present condition for Case 1.  Let $D$ be a minimum connected dominating set of  $G|x$. If $x \in D$, Dominator can play on $G|x$ by first choosing $d_1^*=x$ and  then playing on vertices from $D$. Using Lemma~\ref{lem:1}, this ensures that the game on $G|x$ ends after at most $2 |D| - 1 = 2 \gamma_{\rm c}(G|x) - 1$ moves. If $x \notin D$, then Dominator may play the vertices from $D$ in any legal order and the game on $G|x$ ends after at most $2 \gamma_{\rm c}(G|x) - 1$ moves again. As every connected dominating set of $G$ is a connected dominating set of $G|x$, we have $\gamma_{\rm c}(G|x) \le \gamma_{\rm c}(G)$. This, together with the lower bound in Theorem~\ref{thm:c-bound}, proves~\eqref{eq:weak-bounds}.	
	Assuming either $\gamma_{\rm c}(G|x) \le \gamma_{\rm c}(G)-1 $ or $\gamma_{\rm c}(G) \le \ggcd(G)-1$ (or both), the modified inequality chain results in $\ggcd(G|x) \leq 2 \ggcd(G) -3$ as desired.	
\end{description}

From now on, we suppose that $\gamma_{\rm c}(G|x) = \gamma_{\rm c}(G)$ and $\gamma_{\rm c}(G) = \ggcd(G)$. The latter equality implies that for every $P\in \cG$, the set  $\widehat{P}$ of moves is a minimum connected dominating set in $G$. Therefore, Staller may play arbitrary legal moves during the game and it always corresponds to one of her optimal strategies. 
\begin{description}
	\item[Case 2:] There exists an optimal game $P\in \cG$ such that $\widehat{P}_{2k-3}$ contains a vertex $v$ from $N[x]$. 
	
	Consider the following startegy of Dominator on the predominated graph $G|x$. He first plays $v$ and then chooses vertices from $\widehat{P}_{2k-3}$ while it is possible. Since the first move $v$ is from $N[x]$, the continuation of the game on $G|x$ corresponds to the continuation on $G$. According to Dominator's strategy,  after $2(2k-3)-1=4k-7$ moves a superset of $N[\widehat{P}_{2k-3}]$ is dominated. If the game is not over yet,  for the dominated vertices that are outside of $N[\widehat{P}_{2k-3}]$, we can specify a set $C$ of moves such that  $N[\widehat{P}_{2k-3} \cup C]$ gives exactly the set of vertices dominated during the first $4k-7$ moves of the game on $G|x$. Then, starting with Staller's turn, $\ggcd'((\widehat{P}_{2k-3} \cup C) \ccg G)$ further moves are needed to finish the game. By the Connected Game Continuation Principle, it is at most $\ggcd'(\widehat{P}_{2k-3} \ccg G)$.
	
	If $\ggcd(G)=2k$, then $\ggcd'(\widehat{P}_{2k-3} \ccg G)=3$ and the above strategy of Dominator ensures that the game ends in at most $4k-7+3=2\ggcd(G)-4$ moves on $G|x$. If $\ggcd(G)=2k-1$, then $\ggcd'(\widehat{P}_{2k-3} \ccg G)=2$ and the game ends in at most $4k-7+2=2\ggcd(G)-3$ moves on $G|x$. This verifies the upper bound $2\ggcd(G)-3$ for Case 2. 
	
	\item[Case 3:] There exists an optimal game $P\in \cG$ such that  after the moves $d_1, s_1, \dots, d_{k-1}$ on $G$, every legal move dominates only $x$ from $V(G) \setminus N[\widehat{P}_{2k-3}]$.
	
	Let $Y=N(x) \cap N[\widehat{P}_{2k-3}]$, $Z= N(x) \setminus Y$, and $W=V(G) \setminus ( N[\widehat{P}_{2k-3}] \cup \{x\} )$. By definition, $Z \subseteq W$. Under the condition given for this case, $x$ is a cut vertex as there is no edge between $W$ and $N[\widehat{P}_{2k-3}]$. In the optimal game $P$ on $G$,  $Y$ is the set of legal choices for Staller after $d_{k-1}$ is played and then, the only legal move for Dominator is $d_k=x$.
	
	First, suppose that $\ggcd(G)=2k-1$. Then the move $d_k=x$ dominates $Z$ and finishes the game on $G$. Thus, $Z=W$.  Consider the strategy of Dominator on $G|x$ when his first move is  $d_1^*=x$ and then, he plays vertices only from $\widehat{P}\setminus \{d_k,s_{k-1}\}$. We show that this can be done and that the game on $G|x$ finishes after at most $4k-5$ moves. The main observation is that after $d_1^*=x$, no vertex from $Z=W$ is a legal move. Hence Staller must play a vertex $y \in Y$ as $s_1^*$, and then $(\widehat{P}\setminus \{s_{k-1}\}) \cup \{y\}$ remains a connected dominating set in $G$. The set $(\widehat{P}\setminus \{s_{k-1}\}) \cup \{y\}$ is indeed connected as $y \in N[\hat{P}_{2k-3}]$, $xy \in E(G)$, and all neighbors of $s_{k-1} $ except $x$ have a neighbor in $\hat{P}_{2k-3}$. Then, by Lemma~\ref{lem:2}, Dominator can ensure that the game on $G|x$ finishes in at most $2+ 2(2k-3)-1=4k-5=2\ggcd(G)-3$ moves.
	
	If  $\ggcd(G)=2k$, Staller finishes the game $P$ with the move $s_k$. Similarly to the previous argumentation, as $x$ is a cut vertex and $Y$ contains all the legal moves after $d_{k-1}$, we may infer $d_k =x$. It also follows that $s_k \in Z$. Since each legal move finishes the game on $G$ after $d_k=x$, $s_k$ may be replaced in $P$ by each $z \in Z$ that satisfies $N[z]\setminus N[x] \neq \emptyset$. Let Dominator's first move be $d_1^*=x$ on $G|x$. If Staller replies with a vertex $z \in Z$, then $N[z]\setminus N[x] \neq \emptyset$. Therefore, $s_k$ can be replaced by $z$ in the optimal game $P$, and Dominator continues playing vertices from  $\widehat{P}\setminus \{d_k, s_{k}\}$. By Lemma~\ref{lem:2}, this game on $G|x$ finishes in at most $2+ 2(2k-2)-1=4k-3=2\ggcd(G)-3$ moves. In the other case, Staller chooses a vertex $s_1^*=y$ from $Y$. We may observe again that $s_{k-1}$ can be replaced by $y$ in $P$. That is, if  Dominator always selects a legal move from $\widehat{P}\setminus \{d_k, s_{k-1}\}$ in the continuation, the game finishes after at most $2+ 2(2k-2)-1=4k-3=2\ggcd(G)-3$ moves on $G|x$.
	
	\item[Case 4:] There exists an optimal game $P\in \cG$ such that  after the moves $d_1, s_1, \dots,$ $d_{k-1}, s_{k-1}$ on $G$, every legal move dominates only $x$ from $V(G) \setminus N[\widehat{P}_{2k-2}]$.
	
	Observe first that Case 4 cannot arise if $d_k$ is the last move in the game $P$ on $G$. Indeed, in this case already $\widehat{P} \setminus \{d_k\}$ would be a connected dominating set in $G|x$, contradicting the equality $\gamma_{\rm c}(G|x) = \gamma_{\rm c}(G)$ that must be true under the exclusion of Case 1.
	
	Suppose now that $\ggcd(G)=2k$ and after $s_{k-1}$, every legal move dominates only $x$. As follows, $x$ is a cut vertex, and after $d_k$ is played  the only legal choice for Staller is $s_k=x$ that finishes the game. In this case, every $y \in N(x) \cap N[\widehat{P}_{2k-2}]$ is an optimal move for Dominator and therefore, $d_{k}$ can be replaced by $y$ in $P$. Similarly to the previous case, Dominator starts with the move $d_1^*=x$ on $G|x$. In the next turn, Staller must play a vertex $y$ from  $N(x) \cap N[\widehat{P}_{2k-2}]$. From this point, Dominator always plays a vertex from $\widehat{P}_{2k-2}$ that ensures that the game finishes after at most $2+2(2k-2)-1=4k-3= 2\ggcd(G)-3$ moves on $G|x$.
	
	\item[Case 5:] There exists an optimal game $P\in \cG$ such that, at a stage of the game,  every optimal move of the next player dominates only $x$ from the set of vertices which have not been dominated so far.
	
	Since Case 2 is excluded, it cannot happen during the first $2k-3$ moves of $P$. Since Case 3 is excluded and every legal choice is optimal for Staller, this move cannot be $s_{k-1}$. Further, this move cannot be the last move of the game on $G$, because then, by deleting the last move from $\widehat{P}$, we would obtain a connected dominating set of $G|x$ that is smaller then $\ggcd(G)$. This contradicts the exclusion of Case 1. Therefore, we may assume that $\ggcd(G)=2k$ and Dominator selects the move $d_k$ when this situation arises. Since Case 4 is excluded, we may suppose that there is a legal, but not optimal, move $u$ such that $N(u) \setminus N[\widehat{P}_{2k-2}]$ contains a vertex different from $x$. In fact, by Corollary~\ref{cor:adv}, the non-optimal move $u$ cannot dominate $x$. Now, consider the optimal game $P$ on $G$. If the last move $s_k$ is different from $x$, then $\widehat{P}\setminus \{d_k\}$ is a connected dominating set in $G|x$ that contradicts the equality $\gamma_{\rm c}(G|x)= \ggcd(G)$. Note that in this case $\widehat{P}\setminus \{d_k\}$ is connected because it contains the first $2k - 2$ moves of the game, and since $s_k$ is different from $x$ it must be connected to one of the previous moves. Thus, $s_k=x$ and this move dominates the set $Z=N(x) \setminus N[\widehat{P}_{2k-2}]$. After the move $d_k$, the vertex $u$ remains playable. Since Staller cannot delay the end of the game, $u$ dominates the entire set $Z$. It gives a contradiction again, as $u$ would be an optimal move for Dominator. Indeed, if Dominator selects  $u$ as his $k^{\rm th}$ move, then Staller is forced to finish the game on $G$ by dominating $x$ with her $k^{\rm th}$ move. We conclude that the condition of Case 5 cannot be satisfied if each of Case 1--4 is excluded.
	
	\item[Case 6:] None of the previous conditions is true.
	
	In this case, Dominator plays on $G|x$ by following an optimal startegy for $G$ with the restriction that he never plays a vertex which dominates only $x$. Since Case 5 is excluded, Dominator always has an optimal move on $G$ that satisfies this condition. By the same reason, Staller always has a legal move on $G|x$ that is also optimal on $G$. In this way, Dominator's strategy ensures that the game on $G|x$ finishes in exactly $\ggcd(G)$ moves and consequently, $\ggcd(G|x) = \ggcd(G)$. Using the condition  $\ggcd(G)\ge 3$, we conclude the desired result $\ggcd(G|x) \le 2\ggcd(G)-3$.  \qed
\end{description}

\subsection{Sharpness}

Sharpness of the bound of Theorem~\ref{thm:G|x_upper} follows from the following result. 

\begin{proposition}
	\label{prop:predom+(n-3)}
	For every $n \geq 3$ there exists a graph $G$ with a vertex $x$ such that $$ \ggcd(G) = n \quad \text{and} \quad \ggcd(G|x) = 2 n - 3.$$
\end{proposition}

\proof
	Consider a graph $H_n'$ defined as follows. Take the graph $H_{n}$ defined in Section~\ref{sec:prelim} and remove vertices $x_1$ and $y_1$. Set $U = \{ u_0, \ldots, u_{n+1} \}$. Since $\gamma_c(H_n') = n$ and Dominator can play according to strategy \fast\ as described in Section~\ref{sec:prelim}, we have $\ggcd(H_n') = n$. In the following we prove that $\ggcd(H_n'|u_1) = 2n-3 = 2 \ggcd(H_n') - 3$.
	
	If Dominator does not start the game on $N[u_1]$, then the vertex $u_3$ will be played before $u_2$. But this means that after $u_3$ is played, playing $u_2$ is not legal, thus $u_0$ can never be dominated. Hence, Dominator must start the game on $N[u_1]$ to ensure that the game finishes in a finite number of moves. We consider all three possible start moves of Dominator.
	
	\begin{description}
		\item[Case 1:] Dominator starts on $d_1 = u_2$.\\
		Staller follows the strategy \slow\ (from Section~\ref{sec:prelim}) and plays $u_1$ only if she is forced to do so in the last move of the game. If she is indeed forced to finish the game on $u_1$, then due to strategy \slow\ at least $n-2$ moves are played outside of $U$ and altogether at least $(n-1) + (n-2) + 1 = 2 n - 2$ moves are played. Otherwise, Dominator playes $u_1$ during the game. After this move, Staller can only reply on a vertex from $U$, so she ensures only $n-3$ moves outside $U$. Thus the number of moves in this case is at least $n + (n-3) = 2n-3$.
		
		\item[Case 2:] Dominator starts on $d_1 = u_1$.\\
		Staller is forced to reply on $s_1 = u_2$. After this move, she follows strategy \slow\ to ensure that at least $n-3$ moves are played outside of $U$. Altogether, in this case at least $n + (n-3) = 2n-3$ moves are needed to finish the game on $H_n'|u_1$.
		
		\item[Case 3:] Dominator starts on $d_1 = u_0$.\\
		Since $N[u_0] \subseteq N[u_1]$, it follows from Corollary~\ref{cor:adv} that the number of moves in this case is at least the number of moves in the game where Dominator starts on $u_1$. Thus by Case 2, the number of moves is at least $2n - 3$. 
	\end{description}
	
	We see that independently of Dominator's first move, Staller can ensure that at least $2n-3$ moves are played. Hence $\ggcd(H_n'|x) = 2n-3$ by the already proved upper bound and since $\ggcd(H_n') = n$.
\qed

The graphs $H_n'$ from Proposition~\ref{prop:predom+(n-3)} contain cut-vertices, thus we wonder whether the upper bound of Theorem~\ref{thm:G|x_upper} is also sharp on $2$-connected graphs. While we do not have an answer to this question, we were nevertheless able to see that the difference $\ggcd(G|x) - \ggcd(G)$ can be arbitrarily large also for $2$-connected graphs $G$. For this sake let $C_{1,3}$ be the graph from Fig.~\ref{fig:C_1,3}. Using a computer, we obtain $$\ggcd(C_{1,3}|w) = 16 > 14 =  \ggcd(C_{1,3})\,.$$

\begin{figure}[!ht]
	\begin{center}
		\begin{tikzpicture}[thick,scale=0.1]
		
		\pgfmathtruncatemacro{\h}{10}
		\pgfmathtruncatemacro{\d}{10}
		\pgfmathtruncatemacro{\H}{8}
		\pgfmathtruncatemacro{\D}{4}
		\pgfmathtruncatemacro{\k}{1}
		\pgfmathtruncatemacro{\l}{3}
		
		\node (s) at (0,0) {};
		
		\foreach \x in {0,1,...,\k}
		\node (a\x) at (\x*\d,\h) {};
		\node[label=above: $w$] (x) at (\k*\d + \d, \h) {};
		
		\foreach \x in {0,1,...,\k}
		\node (b\x) at (\x*\d,-\h) {};
		\node (y) at (\k*\d + \d, -\h) {};
		
		\foreach \x in {0,1,...,\l}
		\node (c\x) at (\k*\d + 2*\d + \x*\d,\h) {};
		
		\foreach \x in {0,1,...,\l}
		\node (d\x) at (\k*\d + 2*\d + \x*\d,-\h) {};
		
		\node (t) at (\k*\d + 2*\d + \l*\d,0) {};
		
		\foreach \x [remember=\x as \lastx (initially 0)] in {1,...,\k}
		\path (a\x) edge (a\lastx);
		
		\foreach \x [remember=\x as \lastx (initially 0)] in {1,...,\k}
		\path (b\x) edge (b\lastx);
		
		\foreach \x [remember=\x as \lastx (initially 0)] in {1,...,\l}
		\path (c\x) edge (c\lastx);
		
		\foreach \x [remember=\x as \lastx (initially 0)] in {1,...,\l}
		\path (d\x) edge (d\lastx);
		
		\path (s) edge (a0);
		\path (s) edge (b0);
		\path (x) edge (a\k);
		\path (x) edge (c0);
		\path (y) edge (b\k);
		\path (y) edge (d0);
		\path (t) edge (c\l);
		\path (t) edge (d\l);
		
		\foreach \x in {1,...,\k}
		\node (a'\x) at (\x*\d-\D,\h+\H) {};
		\foreach \x in {1,...,\k}
		\node (a''\x) at (\x*\d,\h+\H) {};
		\foreach \x in {1,...,\k}
		\node (b'\x) at (\x*\d-\D,-\h-\H) {};
		\foreach \x in {1,...,\k}
		\node (b''\x) at (\x*\d,-\h-\H) {};
		
		\foreach \x in {1,...,\l}
		\node (c'\x) at (\k*\d + 2*\d + \x*\d-\D,\h+\H) {};
		\foreach \x in {1,...,\l}
		\node (c''\x) at (\k*\d + 2*\d + \x*\d,\h+\H) {};
		\foreach \x in {1,...,\l}
		\node (d'\x) at (\k*\d + 2*\d + \x*\d-\D,-\h-\H) {};
		\foreach \x in {1,...,\l}
		\node (d''\x) at (\k*\d + 2*\d + \x*\d,-\h-\H) {};
		
		\foreach \a in {a,b}
		\foreach \x [remember=\x as \lastx (initially 0)] in {1,...,\k}
		\draw (\a\x) -- (\a''\x) -- (\a\lastx) -- (\a'\x) -- (\a''\x);
		
		\foreach \a in {c,d}
		\foreach \x [remember=\x as \lastx (initially 0)] in {1,...,\l}
		\draw (\a\x) -- (\a''\x) -- (\a\lastx) -- (\a'\x) -- (\a''\x);

		\end{tikzpicture}
		\caption{The graph $C_{1,3}$.}
		\label{fig:C_1,3}
	\end{center}
\end{figure}
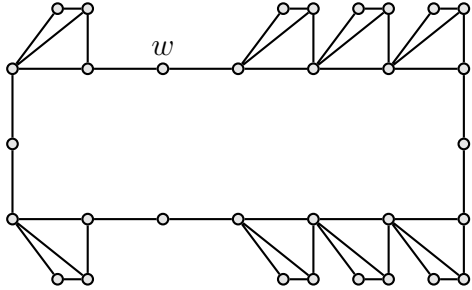

Let now $C_{k, \ell}$ be a graph obtained from two copies $A$ and $B$ of $H_{k+1}$ and two copies $C$ and $D$ of $H_{\ell+1}$ by identifying several vertices. Denote vertices of $A$ by $a_i, a_i', a_i''$ for $u_{k+2-i}, x_{k+2-i}, y_{k+2-i}$ from the definition of the graph $H_{k+1}$, respectively. Similar notation is used for vertices of $B$, $C$ and $D$. To obtain the graph $C_{k, \ell}$, identify the following pairs of vertices: $a_0$ and $b_0$, $c_{\ell+2}$ and $d_{\ell+2}$, $a_{k+2}$ and $c_0$, $b_{k+2}$ and $d_0$, and label them $s$, $t$, $w$, $y$, respectively. Let $C_0=\{s\} \cup \{a_1, \ldots, a_{k+1}\} \cup \{w\} \cup \{b_1, \ldots, b_{k+1}\} \cup \{y\} \cup \{c_1, \ldots, c_{\ell+1}\} \cup \{d_1, \ldots, d_{\ell+1}\} \cup \{t\}$. See Fig.~\ref{fig:C_k,l} for the graph $C_{3,6}$ and note that $C_0$ is the set of the vertices of the inner long cycle.

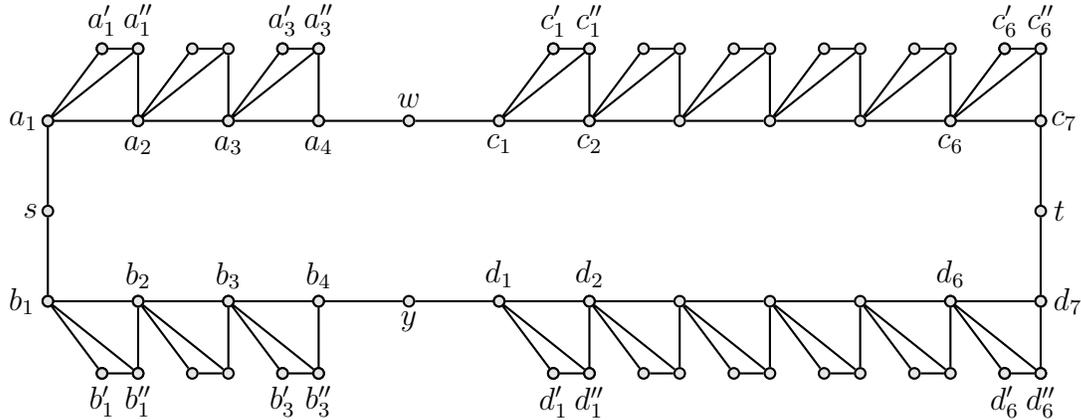
\begin{figure}[!hb]
	\begin{center}
		\begin{tikzpicture}[thick,scale=0.12]
		
		\pgfmathtruncatemacro{\h}{10}
		\pgfmathtruncatemacro{\d}{10}
		\pgfmathtruncatemacro{\H}{8}
		\pgfmathtruncatemacro{\D}{4}
		\pgfmathtruncatemacro{\k}{3}
		\pgfmathtruncatemacro{\l}{6}
		
		\node[label=left: $s$] (s) at (0,0) {};
		
		\foreach \x in {0,1,...,\k}
		\node (a\x) at (\x*\d,\h) {};
		
		\node[label=left:$a_1$] (A1) at (0,\h) {};
		\node[label=below:$a_2$] (A2) at (\d,\h) {};
		\node[label=above:$a'_1$] (A'1) at (\d-\D,\h+\H) {};
		\node[label=above:$a''_1$] (A''1) at (\d,\h+\H) {};
		
		\node[label=below:$a_3$] (Ak) at (\k*\d-\d,\h) {};
		\node[label=below:$a_{4}$] (Ak+1) at (\k*\d,\h) {};
		\node[label=above:$a'_3$] (A'k) at (\k*\d-\D,\h+\H) {};
		\node[label=above:$a''_3$] (A''k) at (\k*\d,\h+\H) {};
		
		\node[label=above: $w$] (x) at (\k*\d + \d, \h) {};
		
		\foreach \x in {0,1,...,\k}
		\node (b\x) at (\x*\d,-\h) {};
		
		\node[label=left:$b_1$] (B1) at (0,-\h) {};
		\node[label=above:$b_2$] (B2) at (\d,-\h) {};
		\node[label=below:$b'_1$] (B'1) at (\d-\D,-\h-\H) {};
		\node[label=below:$b''_1$] (B''1) at (\d,-\h-\H) {};
		
		\node[label=above:$b_3$] (Bk) at (\k*\d-\d,-\h) {};
		\node[label=above:$b_{4}$] (Bk+1) at (\k*\d,-\h) {};
		\node[label=below:$b'_3$] (B'k) at (\k*\d-\D,-\h-\H) {};
		\node[label=below:$b''_3$] (B''k) at (\k*\d,-\h-\H) {};
		
		\node[label=below: $y$] (y) at (\k*\d + \d, -\h) {};
		
		\foreach \x in {0,1,...,\l}
		\node (c\x) at (\k*\d + 2*\d + \x*\d,\h) {};
		
		\node[label=below:$c_1$] (C1) at (\k*\d + 2*\d,\h) {};
		\node[label=below:$c_2$] (C2) at (\k*\d + 2*\d + \d,\h) {};
		\node[label=above:$c'_1$] (C'1) at (\k*\d + 2*\d + \d-\D,\h+\H) {};
		\node[label=above:$c''_1$] (C''1) at (\k*\d + 2*\d + \d,\h+\H) {};
		
		\node[label=below:$c_6$] (Cl) at (\k*\d + 2*\d + \l*\d-\d,\h) {};
		\node[label=right:$c_{7}$] (Cl+1) at (\k*\d + 2*\d + \l*\d,\h) {};
		\node[label=above:$c'_6$] (C'l) at (\k*\d + 2*\d + \l*\d-\D,\h+\H) {};
		\node[label=above:$c''_6$] (C''l) at (\k*\d + 2*\d + \l*\d,\h+\H) {};
		
		\foreach \x in {0,1,...,\l}
		\node (d\x) at (\k*\d + 2*\d + \x*\d,-\h) {};
		
		\node[label=above:$d_1$] (D1) at (\k*\d + 2*\d,-\h) {};
		\node[label=above:$d_2$] (D2) at (\k*\d + 2*\d + \d,-\h) {};
		\node[label=below:$d'_1$] (D'1) at (\k*\d + 2*\d + \d-\D,-\h-\H) {};
		\node[label=below:$d''_1$] (D''1) at (\k*\d + 2*\d + \d,-\h-\H) {};
		
		\node[label=above:$d_6$] (Dl) at (\k*\d + 2*\d + \l*\d-\d,-\h) {};
		\node[label=right:$d_{7}$] (Dl+1) at (\k*\d + 2*\d + \l*\d,-\h) {};
		\node[label=below:$d'_6$] (D'l) at (\k*\d + 2*\d + \l*\d-\D,-\h-\H) {};
		\node[label=below:$d''_6$] (D''l) at (\k*\d + 2*\d + \l*\d,-\h-\H) {};
		
		\node[label=right: $t$] (t) at (\k*\d + 2*\d + \l*\d,0) {};
		
		\foreach \x [remember=\x as \lastx (initially 0)] in {1,...,\k}
		\path (a\x) edge (a\lastx);
		
		\foreach \x [remember=\x as \lastx (initially 0)] in {1,...,\k}
		\path (b\x) edge (b\lastx);
		
		\foreach \x [remember=\x as \lastx (initially 0)] in {1,...,\l}
		\path (c\x) edge (c\lastx);
		
		\foreach \x [remember=\x as \lastx (initially 0)] in {1,...,\l}
		\path (d\x) edge (d\lastx);
		
		\path (s) edge (a0);
		\path (s) edge (b0);
		\path (x) edge (a\k);
		\path (x) edge (c0);
		\path (y) edge (b\k);
		\path (y) edge (d0);
		\path (t) edge (c\l);
		\path (t) edge (d\l);
		
		\foreach \x in {1,...,\k}
		\node (a'\x) at (\x*\d-\D,\h+\H) {};
		\foreach \x in {1,...,\k}
		\node (a''\x) at (\x*\d,\h+\H) {};
		\foreach \x in {1,...,\k}
		\node (b'\x) at (\x*\d-\D,-\h-\H) {};
		\foreach \x in {1,...,\k}
		\node (b''\x) at (\x*\d,-\h-\H) {};
		
		\foreach \x in {1,...,\l}
		\node (c'\x) at (\k*\d + 2*\d + \x*\d-\D,\h+\H) {};
		\foreach \x in {1,...,\l}
		\node (c''\x) at (\k*\d + 2*\d + \x*\d,\h+\H) {};
		\foreach \x in {1,...,\l}
		\node (d'\x) at (\k*\d + 2*\d + \x*\d-\D,-\h-\H) {};
		\foreach \x in {1,...,\l}
		\node (d''\x) at (\k*\d + 2*\d + \x*\d,-\h-\H) {};
		
		\foreach \a in {a,b}
		\foreach \x [remember=\x as \lastx (initially 0)] in {1,...,\k}
		\draw (\a\x) -- (\a''\x) -- (\a\lastx) -- (\a'\x) -- (\a''\x);
		
		\foreach \a in {c,d}
		\foreach \x [remember=\x as \lastx (initially 0)] in {1,...,\l}
		\draw (\a\x) -- (\a''\x) -- (\a\lastx) -- (\a'\x) -- (\a''\x);

		\end{tikzpicture}
		\caption{The graph $C_{3,6}$.}
		\label{fig:C_k,l}
	\end{center}
\end{figure}

If $k$ and $\ell$ increase, then the following proposition proves that the difference $\ggcd(G|x) - \ggcd(G)$ can be arbitrarily large even for $2$-connected graphs $G$.

\begin{proposition}
	\label{prop:c_k,l}
	If $k, \ell \geq 1$ and $\ell \leq k+1$, then $\ggcd(C_{k, \ell}|w) - \ggcd(C_{k, \ell}) \geq \ell - 2$.
\end{proposition}

\proof
	Since $|C_0| = 2k + 2\ell + 8$, we have $\gamma_c (C_{k, \ell}) = 2k + 2\ell + 6$. Thus $\ggcd(C_{k, \ell}) \geq 2k+2\ell +6$.
	
	Consider the following strategy of Dominator on $C_{k, \ell}$. He starts the game on $s$. We consider the pairs $\{a_i, b_i\}$, $i \in [k+1]$, $\{w, y\}$, and $\{c_i, d_i\}$, $i \in [\ell]$. As soon as Staller plays a vertex from one of those pairs, then Dominator replies on the other vertex from the pair. In particular, Staller's first move will be either $a_1$ or $b_1$, and Dominator's reply will be $b_1$ or $a_1$, respectively. Inductively we can see that all such replies of Dominator are legal and that Staller can only play vertices from $C_0$. Note that by the time when $c_\ell$ and $d_\ell$ are played, $2k + 2\ell + 5$ vertices were played. After that Staller is forced to finish the game with her next move. Thus the game ends in at most $2k + 2\ell + 6$ moves. We conclude that $\ggcd(C_{k, \ell}) = 2k+2\ell +6$.
	
	Now we describe an appropriate strategy for Staller on $C_{k, \ell}|w$, depending on the first move of Dominator.
	
	\begin{description}
		\item[Case 1:] The first move of Dominator is in $\{s\} \cup (V(A) \setminus \{ a_{k+1} \}) \cup B \cup \{y\} \cup D \cup \{t\}$.\\
		Since $w$ is predominated and Dominator's first move is not on $N[w]$, the only way to dominate vertex $c_1$ is by playing vertices $c_{\ell+1}, \ldots, c_2$ in this order. Thus Staller can apply strategy \slow\ on the subgraph $C$. Like this she ensures at least $\ell-1$ moves outside $C_0$. On the other hand, at least $|C_0| - 3$ moves on $C_0$ must be played in order to construct a connected dominating set of the graph. Altogether, at least $(2k+ 2\ell + 8 - 3) + (\ell -1) = \ggcd(C_{k, \ell}) + \ell - 2$ moves are played.
		
		\item[Case 2:] The first move of Dominator is in $C \setminus \{ c_1 \}$.\\
		In order to dominate vertices in $A$, vertices $d_{\ell +1}, \ldots, d_1, y, b_{k+1}, \ldots, b_1$ must be played in this order. Thus Staller's strategy \slow\ allows her to play at least $(\ell-1) + (k-1)$ moves outside $C_0$. On the other hand, at least $|C_0| - 3$ moves must be played on $C_0$. Altogether this gives at least $(2k+ 2\ell + 8 - 3) + (k + \ell - 2) = \ggcd(C_{k, \ell}) + k + \ell - 3$ moves on the graph.
		
		\item[Case 3:] The first move of Dominator is in $N[w] = \{a_{k+1}, w, c_1\}$.\\
		To finish the game, at least $|C_0|-2$ moves on $C_0$ must be played. Staller's strategy is to follow the strategy \slow\ when Dominator plays on $A$ in the direction from $a_{k+1}$ to $a_1$, and on $B \cup D$ when Dominator plays in the direction from $d_{\ell+1}$ to $b_1$. This ensures that she can force at least $k-1$ moves outside $C_0$. Altogether, at least $(2k+ 2\ell + 8 - 2) + (k -1) = \ggcd(C_{k, \ell}) + k - 1$ moves are played.  
	\end{description}

	It follows from the case analysis that $\ggcd(C_{k, \ell}|w) \geq \ggcd(C_{k, \ell}) + \min\{ \ell -2, k+\ell -3, k-1 \}$. Since $k \geq 1$ and $\ell \leq k+1$, $\min\{ \ell -2, k+\ell -3, k-1 \} = \ell-2$. Hence, $\ggcd(C_{k, \ell}|w) - \ggcd(C_{k, \ell}) \geq \ell - 2$.
\qed

Note that the bound in Case 3 can be improved, but since we are only interested in a lower bound for $\ggcd(C_{k, \ell}|x)$ and not necessarily in the exact value, more detailed arguments are omitted.

\section{Proof of Theorem~\ref{thm:G|x_lower}}
\label{sec:proof-of-second-theorem}

Recall the statement of Theorem~\ref{thm:G|x_lower}: if $G$ is a connected graph with $n(G) \geq 2$ and $x \in V(G)$, then
\begin{equation}
	\label{eq:lower}
	 \ggcd(G|x) \geq \left \lceil \frac12 \ggcd(G) \right \rceil\,,
\end{equation}
and the bound is sharp. Before proving the result we note that a slightly weaker bound 
\begin{equation*}
\ggcd(G|x) \ge \left \lceil \frac{1}{2} (\ggcd(G) - 1) \right \rceil 
\end{equation*}
is again much easier to obtain. Indeed, let an optimal D-game be played on $G|x$ and let $D$ be the set of vertices played by the end of the game. Then $D$ is connected and dominates $V(G) \setminus \{x\}$. If $x'$ is an arbitrary neighbor of $x$, then $D \cup \{x'\}$ is a connected dominating set of $G$. It follows that $\gamma_c(G) \leq \ggcd(G|x) + 1$. Thus we have 
$$\ggcd(G|x) \geq \gamma_c(G) - 1 \geq \frac{1}{2} (\ggcd(G) - 1)\,,$$ 
where the last inequality follows from Theorem~\ref{thm:c-bound}. 

\subsection{Proof of~\eqref{eq:lower}}

We prove the result using imagination strategy. The real game, in which Staller plays optimally, is a connected domination game with Chooser played on the graph $G$.  Dominator imagines a connected domination game on $G|x$, and plays optimally in it. Dominator selects an optimal move $d_i$ in the imagined game and tries copying it to the real game. If $d_i$ is a legal move in the real game, he plays it. Otherwise, he plays an arbitrary legal move in the real game. Staller replies optimally in the real game by playing $s_i$. If $s_i$ is a legal move in the imagined game, then Dominator imagines Staller plays $s_i$ there. Otherwise, he selects an arbitrary legal move $\widetilde{s_i}$ for Staller in the imagined game. In this case, before the next move of Dominator, Chooser plays $\widetilde{s_i}$ in the real game. 

Let $D_R$ and $D_I$ be sets of played vertices in the real and in the imagined game, respectively. Note that $D_R$ includes Chooser's moves. We prove that the described interpretations of moves are legal and that the following property holds: 
\begin{equation}
\label{eq:1}
\text{after every move and its interpretation, we have } N[D_I] \subseteq N[D_R].
\end{equation}

Property~\eqref{eq:1} clearly holds after the first move of Dominator. Suppose it is true after a move $s_{i-1}$ of Staller. Dominator selects his optimal reply $d_i$ in the imagined game. In particular, this means that $d_i$ has a neighbor in $D_I$, so $d_i \in N[D_I]$. If $d_i$ is legal in the real game, then Dominator copies it there and~\eqref{eq:1} remains valid. Consider now the case that $d_i$ is not legal in the real game. Since~\eqref{eq:1} held before this move, $d_i \in N[D_I] \subseteq N[D_R]$, so $d_i$ is adjacent to a previously played move also in the real game, but it dominates no new vertices. Playing an arbitrary legal move in the real game then maintains~\eqref{eq:1}.

Suppose that~\eqref{eq:1} holds after a move $d_i$ of Dominator. Staller selects her optimal reply $s_i$ in the real game. If $s_i$ is a legal move in the imagined game, Dominator can copy it there and~\eqref{eq:1} remains valid. Otherwise, Dominator imagines Staller played an arbitrary legal move in the imagined game, say $\widetilde{s_i}$. Since~\eqref{eq:1} held before this move, $\widetilde{s_i}$ is connected to an already played vertex also in the real game, thus Chooser can play it. Afterwards, \eqref{eq:1} still holds.

Let $r$ and $i$ denote the number of moves in the real and in the imagined game, respectively, where $r$ does not include Chooser's moves. Let $c$ denote the number of moves of Chooser in the real game.

Since~\eqref{eq:1} holds throughout the game, thus also at the end, we have $r \leq i+1$. As Staller is playing optimally in the real game, Chooser Lemma~\ref{lema:chooser} gives $\ggcd(G)\leq (r+c)+c$. Since Dominator is playing optimally in the imagined game, we get $\ggcd(G|x) \geq i$. We can also obtain an upper bound for the number of Chooser's moves. He makes zero or one move after each move of Staller in the real game. Thus $c \leq \frac{r}{2}$. We distinguish two cases.
\begin{description}
	\item[Case 1:] $r \leq i$.
	
	We can simplify the bound $r \geq \ggcd(G) - 2 c$ to $2 r \geq \ggcd(G)$, so $r \geq \frac12 \ggcd(G)$. Combining all the obtained inequalities, we get $\frac12 \ggcd(G) \leq r \leq i \leq \ggcd(G|x)$.
	
	\item[Case 2:] $r = i+1$.
	
	After the $i^{\rm th}$ move and its interpretation, the imagined game is over, but the real game is not. Since~\eqref{eq:1} holds at this stage of the game as well, we have $N[D_R] \supseteq N[D_I] \supseteq V(G) \setminus \{x\}$. Since the real game is not finished yet, $x \notin N[D_I]$ and $x \notin N[D_R]$. This means that the vertex $x$ did not get dominated during the games (until the last move in the real game is played). Under this assumption, we can prove that after every move and its interpretation, we actually have the property 
	\begin{equation}
	\label{eq:2}
	D_R = D_I.
	\end{equation}
	
	Property~\eqref{eq:2} clearly holds after Dominator's first move. Suppose it holds after Dominator's move $d_i$. Staller selects her optimal reply $s_i$ in the real game and since we have $D_R = D_I$ and we know that $x$ does not get dominated in the course of the imagined game, $s_i$ is also a legal move in the imagined game and~\eqref{eq:2} remains valid. Suppose now that~\eqref{eq:2} holds after Staller's move $s_{i-1}$. Dominator chooses his optimal reply $d_i$ in the imagined game, which is also legal in the real game because $D_R = D_I$. Again, \eqref{eq:2} remains valid.
	
	Thus there are no moves of Chooser, that is, $c = 0$. Using the obtained inequalities, we get $\ggcd(G) \leq r = i+1 \leq \ggcd(G|x) + 1$, that is, $\ggcd(G) - 1 \leq \ggcd(G|x)$. If $\ggcd(G) \geq 2$, this implies $\frac12 \ggcd(G) \leq \ggcd(G|x)$. If $\ggcd(G) = 1$, then since $G$ is not $K_1$, we also have $\ggcd(G|x) = 1$, so we again get $\frac12 \ggcd(G) \leq \ggcd(G|x)$.
\end{description}
Knowing that $\frac12 \ggcd(G) \leq \ggcd(G|x)$, the desired bound follows.
\qed

\subsection{Sharpness}

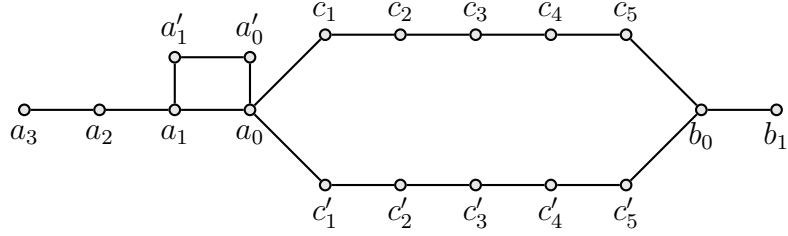
\begin{figure}[!ht]
	\begin{center}
		\begin{tikzpicture}[thick,scale=1]
		
		\node[label=below: $a_0$] (a0) at (0,0) {};
		\node[label=below: $a_1$] (a1) at (-1,0) {};
		\node[label=below: $a_2$] (a2) at (-2,0) {};
		\node[label=below: $a_3$] (a3) at (-3,0) {};
		\node[label=above: $a_0'$] (a0') at (0,0.7) {};
		\node[label=above: $a_1'$] (a1') at (-1,0.7) {};
		
		\draw (a0) -- (a1) -- (a2) -- (a3);
		\draw (a0) -- (a0') -- (a1') -- (a1);
		
		\node[label=above: $c_1$] (c1) at (1,1) {};
		\node[label=above: $c_2$] (c2) at (2,1) {};
		\node[label=above: $c_3$] (c3) at (3,1) {};
		\node[label=above: $c_4$] (c4) at (4,1) {};
		\node[label=above: $c_5$] (c5) at (5,1) {};
		
		\node[label=below: $c'_1$] (d1) at (1,-1) {};
		\node[label=below: $c'_2$] (d2) at (2,-1) {};
		\node[label=below: $c'_3$] (d3) at (3,-1) {};
		\node[label=below: $c'_4$] (d4) at (4,-1) {};
		\node[label=below: $c'_5$] (d5) at (5,-1) {};
		
		\node[label=below: $b_0$] (b0) at (6,0) {};
		\node[label=below: $b_1$] (b1) at (7,0) {};
		
		\draw (a0) -- (c1) -- (c2) -- (c3) -- (c4) -- (c5) -- (b0) -- (b1);
		\draw (a0) -- (d1) -- (d2) -- (d3) -- (d4) -- (d5) -- (b0);
		
		\end{tikzpicture}
		\caption{The graph $D_{5}$.}
		\label{fig:D_5}
	\end{center}
\end{figure}

Examples of graphs with $\ggcd(G|x) = \ggcd(G) - 1$ are already known from~\cite{irsic2019+connected} and include cycles and paths. Another example with the same property are graphs $G$ with $\Delta(G) = n(G) - 2$. If $x \in V(G)$ is the vertex which is not adjacent to the vertex of degree $\Delta(G)$, then $\ggcd(G|x) = 1 = \ggcd(G) - 1$. We were able to obtain an example where the drop is $2$ and generalize it to an infinite family with the same property. Let $D_n$, $n \geq 3$, be a graph obtained in the following way. Take $C_{2n+2}$ on vertices $a_0, c_1, \ldots, c_n, b_0, c'_n, \ldots, c'_1$ and corresponding edges. Add vertices $b_1$, $a_1$, $a_2$, $a_3$, $a_0'$, $a_1'$ and edges $b_0 b_1$, $a_0  a_1$, $a_1 a_2$, $a_2 a_3$, $a_0 a_0'$, $a_0' a_1'$, $a_1' a_1$. See Fig.~\ref{fig:D_5} for $D_5$. Denote $C = \{c_1, \ldots, c_n\}$ and $C' = \{c'_1, \ldots, c'_n\}$. With a case analysis, the following can be obtained.

\begin{proposition}
	\label{prop:-2}
	If $n \geq 3$, then $\ggcd(D_n|c_2) = \ggcd(D_n) - 2$.
\end{proposition}

\begin{figure}[ht!]
	\begin{center}
		\begin{tikzpicture}[scale=1.0,style=thick]
		\tikzstyle{every node}=[draw=none,fill=none]
		
		\def\vr{2pt} 
		
		\begin{scope}[yshift = 0cm, xshift = 0cm]
		\path (0,0) coordinate (v1);
		\path (3,0) coordinate (v2);
		\path (1.5,1.0) coordinate (y1);
		\path (1.5,0.4) coordinate (x1);
		\path (1.5,3.0) coordinate (y2);
		\path (1.5,2.2) coordinate (x2);
		%
		\draw (v1) -- (v2) -- (x1) -- (y1) -- (v1) -- (x1);
		\draw (v1) -- (x2) -- (y2);
		\draw (x2) .. controls (2,2) and (2.5,1.5) .. (v2);
		\draw (y2) .. controls (0.2,2) and (0.2,1) .. (v1);
		\draw (v1)  [fill=white] circle (\vr);
		\draw (v2)  [fill=white] circle (\vr);
		\draw (y1)  [fill=white] circle (\vr);
		\draw (x1)  [fill=white] circle (\vr);
		\draw (y2)  [fill=white] circle (\vr);
		\draw (x2)  [fill=white] circle (\vr);
		\draw[right] (y1)++(0,0.25)  node {$y_1^{(1)}$};
		\draw[right] (x1)++(0,0.25) node {$x_1^{(1)}$};
		\draw[right] (y2)++(0,0.25)  node {$y_r^{(1)}$};
		\draw[right] (x2)++(0,0.25) node {$x_r^{(1)}$};
		\draw[below] (v1) node {$v_1$};
		\draw[above] (y1)++(0,0.3) node {$\vdots$};
		\end{scope}
		
		\begin{scope}[yshift = 0cm, xshift = 3cm]
		\path (0,0) coordinate (v1);
		\path (3,0) coordinate (v2);
		\path (1.5,1.0) coordinate (y1);
		\path (1.5,0.4) coordinate (x1);
		\path (1.5,3.0) coordinate (y2);
		\path (1.5,2.2) coordinate (x2);
		%
		\draw (v1) -- (v2) -- (x1) -- (y1) -- (v1) -- (x1);
		\draw (v1) -- (x2) -- (y2);
		\draw (x2) .. controls (2,2) and (2.5,1.5) .. (v2);
		\draw (y2) .. controls (0.2,2) and (0.2,1) .. (v1);
		\draw (v2) -- (3.6,0);
		\draw (v2) -- (3.5,0.2);
		\draw (v2) -- (3.4,0.4);
		\draw (v2) -- (3.3,0.6);
		\draw (v2) -- (3.2,0.8);
		\draw (v1)  [fill=white] circle (\vr);
		\draw (v2)  [fill=white] circle (\vr);
		\draw (y1)  [fill=white] circle (\vr);
		\draw (x1)  [fill=white] circle (\vr);
		\draw (y2)  [fill=white] circle (\vr);
		\draw (x2)  [fill=white] circle (\vr);
		\draw[right] (y1)++(0,0.25)  node {$y_1^{(2)}$};
		\draw[right] (x1)++(0,0.25) node {$x_1^{(2)}$};
		\draw[right] (y2)++(0,0.25)  node {$y_r^{(2)}$};
		\draw[right] (x2)++(0,0.25) node {$x_r^{(2)}$};
		\draw[below] (v1) node {$v_2$};
		\draw[below] (v2) node {$v_3$};
		\draw[above] (y1)++(0,0.3) node {$\vdots$};
		\end{scope}
		
		\begin{scope}[yshift = 0cm, xshift = 10cm]
		\path (0,0) coordinate (v1);
		\path (3,0) coordinate (v2);
		\path (1.5,1.0) coordinate (y1);
		\path (1.5,0.4) coordinate (x1);
		\path (1.5,3.0) coordinate (y2);
		\path (1.5,2.2) coordinate (x2);
		%
		\draw (v1) -- (v2) -- (x1) -- (y1) -- (v1) -- (x1);
		\draw (v1) -- (x2) -- (y2);
		\draw (x2) .. controls (2,2) and (2.5,1.5) .. (v2);
		\draw (y2) .. controls (0.2,2) and (0.2,1) .. (v1);
		\draw (v1) -- (-0.6,0);
		\draw (v1) -- (-0.55,0.2);
		\draw (v1) -- (-0.4,0.6);
		\draw (v1)  [fill=white] circle (\vr);
		\draw (v2)  [fill=white] circle (\vr);
		\draw (y1)  [fill=white] circle (\vr);
		\draw (x1)  [fill=white] circle (\vr);
		\draw (y2)  [fill=white] circle (\vr);
		\draw (x2)  [fill=white] circle (\vr);
		\draw[right] (y1)++(0,0.25)  node {$y_1^{(n)}$};
		\draw[right] (x1)++(0,0.25) node {$x_1^{(n)}$};
		\draw[right] (y2)++(0,0.25)  node {$y_r^{(n)}$};
		\draw[right] (x2)++(0,0.25) node {$x_r^{(n)}$};
		\draw[below] (v1) node {$v_n$};
		\draw[below] (v2) node {$v_{n+1}$};
		\draw[above] (y1)++(0,0.3) node {$\vdots$};
		\draw[left] (v1)++(-1.5,0) node {$\cdots$};
		\draw[left] (v1)++(-1.5,1.5) node {$\cdots$};
		\end{scope}
		
		\path (7,-1.2) coordinate (z);
		\draw (0,0) -- (z) -- (13,0);
		\draw (0,0)  [fill=white] circle (\vr);
		\draw (13,0)  [fill=white] circle (\vr);
		\draw (z)  [fill=white] circle (\vr);
		\draw[below] (z) node {$z$};
		
		\end{tikzpicture}
	\end{center}
	\caption{The graph $G_{n,r}$.}
	\label{fig:C8}
\end{figure}
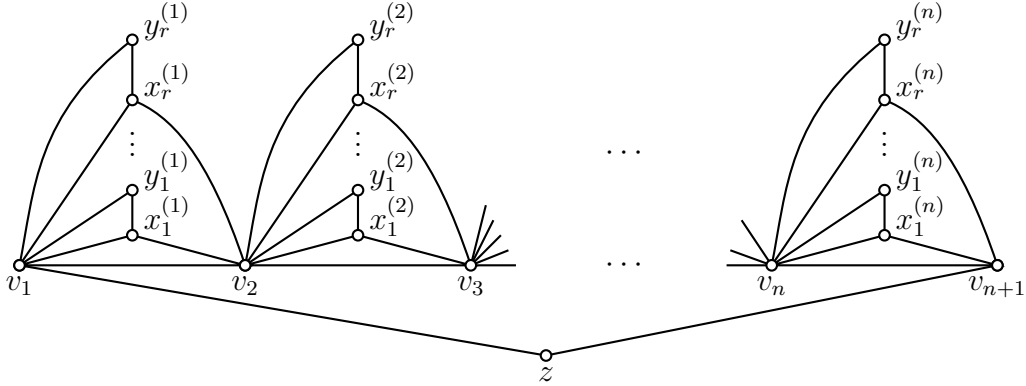

We next present an infinite family that attains the bound from Theorem~\ref{thm:G|x_lower}. Let $G_{n,r}$ be the graph as shown in Figure~\ref{fig:C8}. 
Its vertex set is $\{z\} \cup \{ v_i \; ; i \in [n+1]\} \cup \{ x_j^{(i)} \; ; i \in [n], j \in [r] \} \cup \{ y_j^{(i)} \; ; i \in [n], j \in [r] \}$, and the edges are $z v_1$, $z v_{n+1}$, $v_i v_{i+1}$ for every $i \in [n]$, and the following edges for every $i \in [n]$ and every $j \in [r]$: $v_i x_j^{(i)}$, $v_i y_j^{(i)}$, $x_j^{(i)} y_j^{(i)}$, $x_j^{(i)} v_{i+1}$.

\begin{proposition}
	\label{prop:Gnr}
	If $1 \leq n \leq r$, then 
	$$\ggcd(G_{n,r}) = 2n - 1\quad \text{and} \quad \ggcd(G_{n,r}|v_{n+1}) = n\,.$$ 
\end{proposition}

\proof
Let $V_n = \{v_1, \ldots, v_n, v_{n+1}\}$. Note first that $\gamma_{\rm c}(G_{n,r}) = n$ and that also if $v_{n+1}$ need not to be dominated, we  still need at least $n$ vertices in a connected dominating set for the remaining vertices. It follows that $\ggcd(G_{n,r}|v_{n+1}) \ge n$. Let now Dominator play $v_1$ as his first move on the $G_{n,r}|v_{n+1}$ in the connected domination game. Then the only legal reply for Staller is $v_2$. Dominator then replies with $v_3$. Continuing in this manner we get that the vertices played will be $v_1,  v_2, \ldots, v_n$. This strategy of Dominator yields  $\ggcd(G_{n,r}|v_{n+1}) \le n$.

Consider the D-game played on $G_{n,r}$. As  $\gamma_{\rm c}(G_{n,r}) = n$, we have $\ggcd(G_{n,r}) \le 2n - 1$ by Theorem~\ref{thm:c-bound}. It thus remains to prove that $\ggcd(G_{n,r}) \ge 2n - 1$. For this sake let $T_i$, $i\in [n]$, be the subgraph of $G_{n,r}$ induced by $v_i, v_{i+1}, x_1^{(i)}, \ldots, x_r^{(i)}, y_1^{(i)}, \ldots, y_r^{(i)}$. Suppose that at some point of a D-game, the vertex $v_{i+1}$ has already been played, the vertex $v_i$ was not yet played, and it is Staller's turn. Then we say that Staller has a {\em slow move in $T_i$} if she can play one of the vertices $x_1^{(i)}, \ldots, x_r^{(i)}$. Note that since each of these vertices is adjacent to $v_{i+1}$, such a move $x_\ell^{(i)}$ preserves connectivity and newly dominates only $y_\ell^{(i)}$.  We next describe a strategy of Staller depending on the first move of Dominator. Note that at each stage of the connected game the vertices from $V_n$ that were played so far form an interval $v_i, \ldots, v_{i+k}$, where indices are modulo $n+1$. During the game, at least $n$ moves are played from $V_n$, or if only $n-1$ moves from $V_n$ are played, then at least $r$ further moves are needed.

Suppose first that Dominator plays $z$ as his first move. Then Staller replies with $v_{n+1}$. After that, if Dominator plays $v_1$, Staller plays a slow move in $T_n$, while if Dominator plays $v_n$, Staller plays a slow move in $T_{n-1}$. Proceeding along this way Staller will play at least $n$ moves (which she can since $r \geq n$), hence a total of at least $2n$ moves will be played. 

Suppose second that Dominator plays $v_1$ as his first move. Then Staller replies with $z$. If Dominator then plays $v_2$, Staller replies with $v_{n+1}$. After that, Staller will be able to play at least $n-3$ slow moves, so that she will play in total at least $n-1$ moves. Moreover, after her $(n-3)^{\rm rd}$ slow move, the game is not finished yet, hence at least $2n-1$ moves will be played in total. On the other hand, if Dominator's second move is $v_{n+1}$, then we see as in the previous  case that Staller will play at least $n$ vertices. 

Suppose finally that Dominator plays $v_i$ as his first move, where $i\ge 2$. Then Staller replies with a slow move in $T_{i-1}$. If Dominator then plays $v_{i+1}$, Staller replies with another slow move in $T_{i-1}$. Otherwise, if  Dominator plays $v_{i-1}$, where $i>2$, then Staller replies with a slow move in $T_{i-2}$. Suppose now that at some point of the game Dominator plays $v_1$ and the game is not over yet. Then $z$ is a legal move and Staller replies by playing it. After that Staller's strategy is as in the previous case. This ensures that at least $2n - 1$ moves are played.
\qed

\section{Staller-start game}
\label{sec:S-game}

In this concluding section, we turn our attention to S-game and characterize graphs $G$ and its vertices $x$ for which $\ggcs(G|x) = \infty$.  For this sake, the following concept will be useful.

Let $G$ be a graph and $x\in V(G)$. We say that a player of a connected domination game played on $G$ has an {\em $x$-isolation strategy} if he or she can play such that each vertex at distance $1$ or $2$ from $x$ is dominated before a single vertex from $N_G[x]$ has been played. For instance, if $x$ is the first vertex of a path $P_n$, $n\ge 3$, then in the D-game, Dominator has an $x$-isolation strategy by playing the third vertex of $P_n$ as his first move. Similarly, Staller has an $x$-isolation strategy in the S-game by playing the same vertex as her first move. 

\begin{theorem}
	\label{thm:S-game-finite}
	Let $x$ be a vertex of a connected graph $G$. Then $\ggcs(G|x) = \infty$ if and only if 
	\begin{enumerate}
		\item[(i)] $x$ is a cut-vertex and 
		\item[(ii)]
		$G-x$ contains a component $H$ such that Staller has an $x$-isolation strategy in the S-game played on $G[V(H) \cup \{x\}]$. 
	\end{enumerate}
\end{theorem}

\proof
The result is clearly true for complete graphs $K_1$ and $K_2$, hence assume in the rest that $n(G) \ge 3$. 

Suppose first that $x$ is not a cut-vertex. Let $G'$ be the block of $G$ containing $x$. (It is possible that $G=G'$.) Consider now the S-game played on $G|x$. Since $n(G)$ is at least $3$, the first move of Staller is always possible. When the game continues, eventually a vertex from $G'$ must be played. We claim that all the vertices of $G'$ will be dominated during the game. For this sake let $u\in V(G')$ be a vertex which is not yet dominated by the vertices played so far. Let $w$ be an arbitrary vertex of $G'$ that has already been played, and note that $wu\notin E(G)$. (It is possible that $w=x$.) Let $P$ and $P'$ be internally disjoint $u,w$-paths in $G'$. These paths exist as $G'$ is $2$-connected. Note also that both paths are of length at least $2$ because $wu\notin E(G)$. At least one of $P$ and $P'$, say $P$, does not contain $x$ as an internal vertex. Let $w'$ be the last vertex on $P$ that has already been dominated during the game played so far. This means that $w'$ has been dominated but not yet played. It follows that $w'$ is a legal move at this moment of the game. By induction, all the vertices of $G'$ will eventually be dominated by the end of the game. Since the only difference between the game played on $G|x$ and on $G$ is the predomination of $x$, this also means that the game will end in a finite number of moves on $G$. We conclude that $\ggcs(G|x) = \infty$ is possible only if $x$ is a cut-vertex. 

Let now $x$ be a cut-vertex. Suppose first that (ii) holds, that is, $G-x$ contains a component $H$ such that Staller has an $x$-isolation strategy in the S-game played on $G[V(H) \cup \{x\}]$. Then Staller plays the first move in $V(H)$ such that her above strategy on $G[V(H) \cup \{x\}]$ is maintained. After achieving her goal, the connectivity condition implies that no vertex from $N_G[x]$ will be played in the rest of the game. Since $x$ is a cut-vertex, it has a neighbor $y\notin V(H)$. But then $y$ cannot be dominated during the game, hence  $\ggcs(G|x) = \infty$. Conversely, suppose that the condition (ii) is not fulfilled. If Staller selects $s_1' = x$, we clearly have $\ggcs(G|x) < \infty$. Otherwise, no matter in which component of $G-x$ Staller starts the game, say in component $H$, Staller cannot prevent Dominator to achieve his goal that some vertex from $N_H(x)$ is played during the game. But then one of the players will eventually need to play $x$. After that the rest of the game has no effect of $x$ being predominated, so it will end after a finite number of moves. 
\qed

We next give a necessary condition for Staller to be able to achieve condition (ii) of Theorem~\ref{thm:S-game-finite}. Using it we are then able to give an explicit description of trees $T$ and its vertices $x$ for which $\ggcs(T|x) = \infty$ holds.

\begin{proposition}
	\label{prp:connected}
	Let $x$ be a cut-vertex of a connected graph $G$ and let $H$ be a component of $G-x$. If $G[V(H)-N_H(x)]$ is disconnected, then Staller does not have an $x$-isolation strategy in the S-game played on $G[V(H)\cup \{x\}]$. 
\end{proposition}

\proof
Let an S-game be played on $G[V(H)\cup \{x\}]$. Suppose that $H' = G[V(H)-N_H(x)]$ is disconnected. Suppose that Staller starts the game by playing  a vertex from a component $H_1$ of $H'$. Let $H_2$ be another component of $H'$. Then, in order that the vertices from $H_2$ become dominated, at least one vertex from $N_H(x)$ will have to be played. But this means that Staller cannot achieve an $x$-isolation strategy.
\qed

\begin{proposition}
	\label{prop:trees}
	Let $x$ be a vertex of a tree $T$. Then $\ggcs(T|x) = \infty$ if and only if  $x$ is not a leaf and has a neighbor of degree $2$.
\end{proposition}

\proof
Suppose that $x$ is not a leaf and that it has a neighbor $x'$ of degree $2$. Let $x''$ be the other neighbor of $x'$.  In an S-game, Staller can start by playing $x''$ which provides an $x$-isolation strategy in the component of $T-x$ which includes $x'$. Thus by Theorem~\ref{thm:S-game-finite} we have $\ggcs(T|x) = \infty$.

Conversely, suppose that $\ggcs(T|x) = \infty$. Then $x$ is not a leaf by Theorem~\ref{thm:S-game-finite}(i). If $x$ would not have a neighbor of degree $2$, then Proposition~\ref{prp:connected} together with Theorem~\ref{thm:S-game-finite} would yield that $\ggcs(T|x)$ is finite. 
\qed

\end{document}